\documentclass[]{imsart}

\RequirePackage{amsthm,amsmath,amsfonts,amssymb}
\RequirePackage[numbers,sort&compress]{natbib}
\RequirePackage[colorlinks,citecolor=blue,urlcolor=blue]{hyperref}
\RequirePackage{graphicx}

\startlocaldefs
\theoremstyle{plain}

\newtheorem{theorem}{Theorem}[section]
\newtheorem{lemma}[theorem]{Lemma}
\newtheorem{corollary}{Corollary}[section]
\newtheorem{remark}{Remark}[section]

\theoremstyle{remark}
\newtheorem{definition}[theorem]{Definition}
\newtheorem{example}{Example}

\newcommand{\R}{\mathbb{R}}
\newcommand{\Nat}{\mathbb{N}}

\endlocaldefs

\begin{document}

\begin{frontmatter}
\title{Entropy bounds for the absolute convex hull of tensors}
%\title{A sample article title with some additional note\thanksref{t1}}
\runtitle{Entropy of tensors}
%\thankstext{T1}{A sample additional note to the title.}

\begin{aug}
%%%%%%%%%%%%%%%%%%%%%%%%%%%%%%%%%%%%%%%%%%%%%%%
%% Only one address is permitted per author. %%
%% Only division, organization and e-mail is %%
%% included in the address.                  %%
%% Additional information can be included in %%
\author[A]{\fnms{Sara}~\snm{van de Geer}\ead[label=e1]{geer@stat.math.ethz.ch}},
%\author[B]{\fnms{Second}~\snm{Author}\ead[label=e2]{second@somewhere.com}\orcid{0000-0000-0000-0000}}
%\and
%\author[B]{\fnms{Third}~\snm{Author}\ead[label=e3]{third@somewhere.com}}
%%%%%%%%%%%%%%%%%%%%%%%%%%%%%%%%%%%%%%%%%%%%%%
%% Addresses                                %%
%%%%%%%%%%%%%%%%%%%%%%%%%%%%%%%%%%%%%%%%%%%%%%
\address[A]{Department of Mathematics,
ETH Z\"urich
\printead[presep={,\ }]{e1}}

%\address[B]{Department,
%University or Company Name\printead[presep={,\ }]{e2,e3}}
\end{aug}

\centerline{February 12, 2024} 
\begin{abstract}
We derive entropy bounds for the absolute convex hull of vectors $X= (x_1 , \ldots , x_p)\in \R^{n \times p} $ in $\R^n$ and apply
this to the case where $X$ is the $d$-fold tensor matrix
$$X = \underbrace{\Psi \otimes \cdots \otimes \Psi}_{d \ {\rm times} }\in \R^{m^d \times r^d },$$ 
with a given
$\Psi = ( \psi_1 , \ldots , \psi_r ) \in \R^{m \times r} $, normalized to that $ \| \psi_j \|_2 \le 1$ for all $j \in \{1 , \ldots , r\}$.
 For $\epsilon >0$ we let 
${\cal V} \subset \R^m$ be the  linear space with smallest dimension 
$M ( \epsilon , \Psi)$ such that $ \max_{1 \le j \le r } \min_{v \in {\cal V} } \| \psi_j -v \|_2 \le \epsilon$.
We call $M( \epsilon , \psi)$ the $\epsilon$-approximation of $\Psi$ and assume it is - up to log terms - polynomial in
$\epsilon$. We show that the entropy of  the absolute convex hull of the $d$-fold tensor matrix $X$
is up to log-terms of the same order as the entropy for the case $d=1$.
The results are generalized to absolute convex hulls of tensors of functions
in $L_2 (\mu)$ where $\mu$ is Lebesgue measure on $[0,1]$.  As an application we consider the space of
functions on $[0,1]^d$ with bounded $q$-th order Vitali total variation for a given $q \in \Nat$. As a by-product, we construct
an orthonormal, piecewise polynomial, wavelet dictionary for functions that are well-approximated by piecewise polynomials.
\end{abstract}

\begin{keyword}[class=MSC]
\kwd[Primary ]{41A10}
\kwd{41A25}
\kwd{94A12}
\kwd[; secondary ]{62G08}
\end{keyword}

\begin{keyword}
\kwd{absolute convex hull}
\kwd{entropy}
\kwd{linear approximation}
\kwd{multi-resolution}
\kwd{piecewise polynomials}
\kwd{tensor}
\kwd{Vitali total variation}
\end{keyword}

\end{frontmatter}
%%%%%%%%%%%%%%%%%%%%%%%%%%%%%%%%%%%%%%%%%%%%%%
%% Please use \tableofcontents for articles %%
%% with 50 pages and more                   %%
%%%%%%%%%%%%%%%%%%%%%%%%%%%%%%%%%%%%%%%%%%%%%%
%\tableofcontents

\section{Introduction}\label{introduction.section}

Let $\Psi= ( \psi_1 , \ldots , \psi_m ) \in \R^{m \times r } $ be an $m \times r$ matrix and
\begin{eqnarray*}
X &:=& \underbrace{\Psi \otimes  \cdots \otimes \Psi }_{d \ {\rm times} }\\
%&=&
%( ( \psi_1 \otimes \psi_1 ,  \ldots , \psi_1 \otimes \psi_q)  , \ldots,
%(\psi_m \otimes \psi_1 , \ldots , \psi_m \otimes \psi_m ) ) \in  \R^{m^2 \times q^2}.
\end{eqnarray*}
the $d$-fold tensor matrix of copies of $\Psi$.
The absolute convex hull of the columns of $X$ is
$$ {\cal F}:= {\rm absconv} (X) := \{ f= X \beta : \ \beta \in \R^{r^d},  \ \| \beta \|_1 \le 1 \} ,$$
where $\| \cdot \|_1 $ denotes the $\ell_1$-norm (see Subsection \ref{notation.section} for details
on the notation). 
We derive bounds for the entropy of ${\cal F}\subset \R^{m^d}$ endowed with the Euclidean metric,
for the case where the columns of the matrix $\Psi$ can be well-approximated
by vectors in a lower-dimensional space.
Entropy bounds based on low-dimensional approximations are given
in \cite{van2023lasso}. We present an improvement
in terms of log-factors in Theorem \ref{entropy.theorem} below. The work in the present paper goes into finding the low-dimensional approximations for tensors.
We develop the arguments in detail in the proof of Theorem \ref{tensorcovering.theorem} for the case $d=2$, as 
the case $d$ possibly larger than 2, given in Corollary \ref{dfold.corollary}, is an easy
extension of these arguments. One sees that the entropy bounds depend on the dimension $d$ only  in the log-terms.
We use a multi-resolution type of approach, namely nested sequences
of linear spaces $\{{\cal V}_k\}_{k \in \Nat_0}$ where each ${\cal V}_k$ has dimension proportional
to $2^k$: see Theorem \ref{tensorcovering.theorem}. 

The entropy results can for example be applied in a regression context to obtain rates of convergence
for tensor denoising, improving those in \cite{ortelli2022tensor} called there ``not so slow rates".

We discuss infinite-dimensional settings as well, where instead of vectors, the objects under study
are functions in $L_2(\mu  \times \cdots \times \mu)$ where $\mu$ is a measure on (say) $\R$. 
This is for instance of interest when studying mixtures of densities as explained in the following example (where $d=2$).

\begin{example} \label{mixtures.example} Let $U= (U_1 , U_2)$ and $V= (V_1 , V_2)$ be random variables in $\R^2$. We observe $U$
but $V$ is an unobserved latent variable (or confounder). We assume that $U_1$ and $U_2$ are conditionally
independent given $V=(v_1 , v_2)$ and that their conditional densities, with respect to a $\sigma$-finite dominating measure $\mu$,
are given by $\psi_{v_1}$ and $\psi_{v_2}$
where $\{ \psi_{\theta} : \ \theta  \in \Theta \} $ is a given family of densities and $\Theta \subset \R$.
Let $G$ be the unknown distribution of $V$. Then the density of $U$ with respect to $\mu \times \mu$ is
$$f (u_1 , u_2) =  \int \psi_{v_1} (u_1) \psi_{v_2} (u_2) d G (v_1 , v_2) , \  (u_1 , u_2) \in \R^2.  $$
\end{example}

Another example, which will be studied in detail in Section  \ref{piecewise.section}, is the following.

\begin{example}\label{piecewise.example}
Let 
$$ {\cal F} := \{ f \in L_2 (\mu \times \cdots \times \mu ): \ \| Df \|_1 \le 1 \} $$
where $\mu$ is Lebesgue measure on $[0,1]$, $Df$ is defined as 
$$D f(u) := \prod_{k=1}^d \partial^q f(u)/  ( \partial u_k)^q , \ u= (u_1 , \ldots , u_d) \in [0,1]^d , $$ and
where $\| Df \|_1 := \int |Df| d (\mu  \times \cdots \times \mu)$ is the $q$-th order Vitali total variation of $f$. 
For the case $q=1$ we recover in Subsection \ref{q0=1.section} the entropy bounds 
for the class of $d$-dimensional distribution functions, as given
in \cite{blei2007metric}. For $q>1$ our - up to log terms tight - entropy bounds appear to be new.
\end{example}
The vector version of the situation described in Example \ref{piecewise.example} is studied in \cite{ortelli2022tensor}.
The functions (or vectors) are linear combinations of spline type functions, called the truncated power functions, and
therefore the theory can be applied to obtain rates for certain multivariate  adaptive
regression splines (MARS) as presented in \cite{friedman1991multivariate}.

The formulation in terms of vectors is useful in the context of regression, where functions
are only evaluated at a finite number of design points. The formulation in terms of functions
however is easier to describe, at least in the context of Vitali total variation: the derivatives and
integrals are easier to present than finite differences and sums. 

As a by-product, we construct in Subsection \ref{multiresolution.section} a orthonormal multi-resolution basis for the space $L_2 (\mu)$ 
of functions on $[0,1]$ endowed with Lebesgue measure $\mu$, that consists of piecewise polynomials
of a given degree. 

In this paper we do not expose the application to regression problems: the focus is on entropy bounds.
In the regression context, one may use the low-dimensional approximations directly, avoiding the route via entropy,
see  \cite{van2023lasso}. This theme will be the content of a separate paper, where we will detail Example 
\ref{mixtures.example} on mixtures.

We have chosen to first describe the problem for vectors, and then state the extension
to functions without proof. One could also state the results in a very general context and then 
look at vectors and functions as special cases. We decided against this approach to avoid beginning 
the paper with abstract
statements. Instead we start with a more concrete setting. We hope to have facilitated the reading
in this way.
Nevertheless, it is of  interest to develop the theory in more general terms (possibly
in terms of Banach and Hilbert spaces) as this may lead to important novel findings.

\subsection{Organization of the paper}\label{organization.section}
After the basic notations and definitions in Section \ref{notation.section}, we present approximation numbers
and $\epsilon$-approximations
in Section \ref{approximation.section} and in Theorem \ref{entropy.theorem} the derived bound for the entropy of
absolute convex hulls. Section \ref{approximationtensors.section} has the $\epsilon$-approximations for
tensors in Theorem \ref{tensorcovering.theorem} and as a consequence, by applying Theorem \ref{entropy.theorem},
the entropy bounds for the absolute convex hulls for tensors in Corollary \ref{dfold.corollary}. Section 
\ref{piecewise.section} studies the case where $\Psi$ consist of heaviside functions
or of hinge functions, or more generally of functions of the
truncated power basis. It presents approximation numbers in Lemma \ref{nestedMARS.lemma}
and the resulting entropy bound, following from Theorem \ref{entropy.theorem} , in Corollary \ref{MARS.corollary}.
Section \ref{conclusion.section} concludes and puts the results in perspective.
Section \ref{proofentropy.section} has the proof of Theorem \ref{entropy.theorem} and Section \ref{proof.section}
the proof of Theorem \ref{tensorcovering.theorem}. Finally, Section \ref{proofspiecewise.section} contains the proofs for
Section \ref{piecewise.section}.

\section{Notation and definitions}\label{notation.section}

For a vector $v \in \R^n$, $\| v \|_2 = (\sum_{i=1}^n v_i^2 )^{1/2}$ denotes the Euclidean (or $\ell_2$-) norm,
$\| v \|_1 := \sum_{i=1}^n | v_i| $ is the $\ell_1$-norm, and $\| v \|_{\infty}:= \max_{1 \le i \le n }| v_i |$ is
the $\ell_{\infty}$-norm.  
%For a matrix $A$, we let $\| A \|_2 $ be the Frobenius norm.

We let ${\rm l}_A $ or ${\rm l} \{ A \}$ be the indicator function of the set $A$. For $u \in \R$,
we write $u_+ := \max\{ u , 0 \}$ and we sometimes write ${\rm l }_{\{ u \ge 0 \}} $ as
$u_+^0 $. 

For $n_1$ and $n_2$ in $\Nat_0$, $n_1 < n_2$, we write $[n_1: n_2]$ for the set
$\{ n_1 , n_1+1 , \ldots , n_2 \} $. 
We write the columns of a matrix $\Psi\in \R^{m \times r}$ as $\psi_j$, $j \in [1: r] $ so that $\Psi= (\psi_1 , \ldots , \psi_r) $. For a matrix $X \in \R^{ n \times p }$ the columns are denoted by $x_j$, $j \in [1: p] $.
We use the same notation for the matrix $\Psi= (\psi_1 , \ldots , \psi_r)$ and the
set of its columns $\{\psi_1 , \ldots , \psi_r\}$ and similarly for the matrix $X$. 
We assume throughout the normalization
$$\max_{1 \le j \le r} \| \psi_j\|_2 \le 1.$$
Note also that then the columns of
\begin{eqnarray*}
\Psi \otimes  \Psi =
( ( \psi_1 \otimes \psi_1 , \ldots ,  \psi_1 \otimes \psi_r)  , \ldots,
(\psi_r \otimes \psi_1 , \ldots , \psi_r \otimes \psi_r ) ) \in  \R^{m^2 \times r^2}
\end{eqnarray*}
have $\ell_2$-norm  at most  1 as well. 

For a function $f: \R^d \rightarrow \R$ we write its argument as $(u_1 , \ldots , u_d)$, i.e.,
the function is $(u_1 , \ldots , u_d) \mapsto f(u_1 , \ldots , u_d)$. We do not apply the standard
notation for the argument $(x_1 , \ldots , x_d) \mapsto f(x_1 , \ldots , x_d)$, not to clash with the notation
$(x_1 , \ldots , x_p)$ for the columns of $X$.

For $n \in \Nat$ and a linear space ${\cal V} \subset \R^n $ we write ${\cal V}^{\bot} $ for its orthocomplement, so that ${\cal V} \oplus {\cal V}^{\bot} = \R^n $. For a vector $f \in \R^n$ we let $f^{\cal V}$ be its projection on ${\cal V}$ so that 
$f = f^{\cal V} + f^{{\cal V}^{\bot}} $. We let ${\rm dim} ({\cal V})$ be the dimension of ${\cal V}$.
For a collection of vectors $\{ v_1 , \ldots , v_k\} $ we let ${\rm span} ( \{ v_1 , \ldots , v_k\} )$ the space spanned by these
vectors, i.e., the collection of all linear combinations of the vectors in $\{ v_1 , \ldots , v_k \} $.
These notations have their obvious counter parts for functions on $\R^d$ in $L_2 ({\bf Q})$ where ${\bf Q}$ is a measure
on $\R^d$.

If $g$ and $h_{n,p}$ are functions on $(0, 1]$, with $h_{n,p}$ possibly depending on $n$ and $p$, but
$g$ not, we write
$h_{n,p} (\epsilon) \lesssim g(\epsilon)$ if 
$$\limsup_{\epsilon \downarrow 0} { h_{n,p} (\epsilon) \over  g (\epsilon)}   \le C< \infty ,$$
where $C$ does not depend on $n$ and $p$.
%Moreover we write 
%$h (\epsilon) \stackrel{\log}{\lesssim} g (\epsilon)$ if there is a constant $c\ge 0$ such 
%that
%$$\limsup_{\epsilon \downarrow 0} { h (\epsilon)  \over  \log^c (1/\epsilon) g (\epsilon) }< \infty .$$

\begin{definition} \label{entropy.definition} Let $T= (T,d) $ be a subset of a metric space. For $\epsilon >0$, the $\epsilon$-covering
number $N( \epsilon , T)$ of $T$ is defined as the minimum number of balls with radius $\epsilon$, necessary to cover
$T$. The entropy of $T$ is ${\cal H} (\cdot , T) := \log N (\cdot , T) $. \\
The local $\epsilon$-entropy of $T$ is defined as
$$ {\cal H} ^{\rm loc} ( \epsilon , T ) := \sup_{t_0 \in T} {\cal H}  (\epsilon/4, \{ t \in T:\ d(t , t_0)  \le \epsilon \} ). $$
\end{definition}

In \cite{yang1999information} local entropies are defined
in terms of the so-called capacity of sets. To avoid digressions, we reformulated it in 
Definition \ref{entropy.definition} in terms of
entropies instead of capacities, at the price of a factor $4$.
Then the result in Lemma 3 of \cite{yang1999information} reads
\begin{equation}\label{localentropy.equation}
{\cal H} (\epsilon/4 , T) - {\cal H} (\epsilon, T) \le  {\cal H} ^{\rm loc} ( \epsilon , T ) . 
\end{equation}

\section{$\epsilon$-approximations and entropy of absolute convex hulls}\label{approximation.section}

In this section we examine the entropy of ${\cal F} := {\rm absconv} (X)$ with $X \in \R^{n \times p}$ not necessarily a $d$-fold tensor matrix. One may also start with a symmetric convex set ${\cal F}$ containing the 
origin and let $X$ be the set of its extreme points.
\begin{definition}\label{approximationnumber.definition}
For $n \in \Nat$ and a linear space ${\cal V} \subset \R^n$ we define its approximation number
at $\Psi$  as
$$ \delta ({\cal V} , X) := \max_{1 \le j \le p } \| x_j^{{\cal V}^{\bot} } \|_2 . $$
\end{definition} 

\begin{definition} \label{epsilonapproximation.definition} For $N \in \Nat$, set $\delta_N (X):= \min \{ \delta ( {\cal V}, X): \ {\rm dim} ({\cal V} ) = N \} $. 
For $\epsilon >0$ the $\epsilon$-approximation of $X$ is defined as
$$ M (\epsilon , X) := \min \{ N: \delta_N (X)  \le \epsilon \} . $$

\end{definition}

The following theorem is an improvement in terms of log-factors of Theorem 5 in \cite{van2023lasso}. 
The idea to bound the local entropy, and then invoke a lemma from  \cite{yang1999information} to show that local and
global entropy are of the same order.
We present a proof in Section \ref{proofentropy.section}. 

\begin{theorem}\label{entropy.theorem} Suppose that $\max_{1 \le j \le p} \| x_j \|_2 \le 1$,
that for some constant $V>0$,
$$ N (\epsilon , X) \lesssim \epsilon^{-V}, $$ and moreover
that for some constants $W>0$ and $w \ge 0$ 
$$ M (\epsilon, X) \lesssim \epsilon^{-W} \log^w (1/\epsilon) . $$
%Then for all $f_0 \in \R^n$ and for $ K_{\epsilon} \in \Nat$  
%\begin{eqnarray*}
%{\cal H} ( \epsilon , \{f \in {\cal F} , \| f - f_0 \|_2 \le 2^{K_{\epsilon}}  \epsilon\}  ) 
%\lesssim  \epsilon^{-{2W \over 2+W} } \log^{2w \over 2+W} (1/\epsilon) \log^{W \over 2+W} (1/ \epsilon) 
% \log^{ 2 \over 2+W}  (2^{K_{\epsilon}})
%\end{eqnarray*}
%and 
Then
\begin{eqnarray*}
{\cal H} (\epsilon, {\cal F} ) &\lesssim &
\epsilon^{-{2W \over 2+W} } \log^{2w \over 2+W} (1/\epsilon) \log^{W \over 2+W} (1/ \epsilon) .
 \end{eqnarray*}
\end{theorem}

\begin{remark} Note that since $M(\epsilon,X) \le N(\epsilon, X)$ for all $\epsilon >0$, we can always take $W \le V$ in the above theorem. 
\end{remark}

\begin{remark}
The assumption $N( \epsilon , X) \lesssim \epsilon^{-V}$ implies
${\cal H} ( \epsilon , {\cal F}) \lesssim \epsilon^{-{2V \over 2+V}}$, see \cite{Ball:90}. 
There is no log-term here, but this bound is not tight in general. 
The bound in Theorem \ref{entropy.theorem} above based on $M (\epsilon, X)$ is however tight up to log-terms, as
explained in \cite{van2023lasso}. 

\end{remark}

\section{$\epsilon$-approximations for tensors}\label{approximationtensors.section}

\begin{definition} \label{nested.definition}
Given a starting value $q \in \Nat$,
let ${\cal V}_0 $ be a $q$-dimensional linear subspace of $\R^m$.  
For $k= 1,2 , \ldots$ let ${\cal V}_k \supset {\cal V}_{k-1}$ be a linear space with dimension $q 2^k$.
We call $\{ {\cal V}_k\}_{k \ge 0} $ a nested sequence of linear spaces with starting value $q$.
\end{definition}

\begin{definition} \label{nestedapproximation.definition}
Let $q \in \Nat$, $\gamma >0$, and $W>0$. We call a nested sequence of linear spaces $\{ {\cal V}_k\}_{k \ge 0}$ 
  a nested approximation for $\Psi$ with parameters $(q, \gamma ,W)$ if it has
  starting value $q$ and for all $k \in \Nat $
$$\delta ({\cal V}_k , \Psi) \le \gamma 2^{-k/W} . $$ 
\end{definition}

The proof of the next theorem was inspired by an approach in \cite{bass1988probability}. 
The latter studies small ball estimates for the Brownian sheet, see Section \ref{functionspaces.section}
for a brief discussing of the connection with entropy estimates.

\begin{theorem} \label{tensorcovering.theorem} Suppose that there exists a nested sequence of approximations for $\Psi$, with parameters
$(q,\gamma,W)$.  Then 
$$ M( \epsilon , \Psi \otimes \Psi ) \lesssim \epsilon^{-W}  \log^{2+W\over 2} ( 1/ \epsilon). $$

\end{theorem} 

\begin{corollary} \label{dfold.corollary}
Recall that the  proof of Theorem \ref{tensorcovering.theorem} is given in Section \ref{proof.section}.
There one sees that the  situation is easily extended to the $d$-fold tensor with $d \ge 2$
(see Remark \ref{d>2.remark}).
Let 
$$X:= \underbrace{\Psi \otimes \cdots \otimes \Psi}_{d \ {\rm times}} . $$
Under the conditions of Theorem \ref{tensorcovering.theorem}, $$ M( \epsilon , X) \lesssim  \epsilon^{-W}  \log^{(d-1)(2+W)\over 2} ( 1/\epsilon)  . $$
If we assume in addition that $N( \epsilon , \Psi) \lesssim \epsilon^{-V}, $ for some constant $V$,
then $ N(\epsilon , X) \lesssim \epsilon^{-dV} $, 
and from Theorem \ref{entropy.theorem} we therefore can conclude that
$$ {\cal H} (\epsilon, {\cal F} ) \lesssim \epsilon^{-{2W \over 2+W}} \log^{d-1} (1/\epsilon) \log^{W \over 2+W} (1/\epsilon) .$$
\end{corollary}

\begin{remark} It is not difficult to generalize to an anisotropic case where
$X= \Psi_1 \otimes \Psi_2$ with 
$\Psi_1$ as well as $\Psi_2$ having nested sequences of approximations with different
parameters, say
$(q_{1},\gamma_{1},W_1)$ and 
$(q_{2},\gamma_{2} ,W_2)$ with $W_2 > W_1$. Then
$$ M( \epsilon , \Psi_1 \otimes \Psi_2 ) \lesssim  \epsilon^{-W_2 }\log^{W_2 \over 2} (1/ \epsilon) . $$
We have a log-term less, but hidden in the constant is the term ${1 /(W_2-W_1)} $.
\end{remark}

\section{Extension to function spaces}\label{functionspaces.section}
Let $X \in \R^{n \times p }$ be given.
One may observe that for $\| \beta \|_1 \le 1$ the mapping
$\beta \mapsto X\beta$ is a mapping from the unit ball in the Banach space $(\R^p,  \| \cdot \|_1)$
to the Hilbert space $(\R^n , \| \cdot \|_2 )$. The dual map is: for $y \in \R^n$, $\| y \|_2 \le 1 $,
$ y \mapsto X^Ty$ which is a mapping from the unit ball in the Hilbert space
$( \R^n , \| \cdot \|_2)$ to the Banach space $( \R^p, \| \cdot \|_{\infty})$. 
Such mappings and the entropy of their range are studied in very general terms,
see \cite{pisier1989volume} and for the connection with small ball estimates, see
\cite{kuelbs1993metric} and  \cite{li1999approximation}. Our approach towards entropy bounds based
on low-dimensional approximations is closely related, but different in the sense
that we avoid the Duality Theorem (derived in \cite{artstein2004duality}, \cite{artstein2004convexified})
and the step involving small ball estimates for Gaussian measures. The approach with small ball
estimates yields tight entropy bounds for the context we consider here, provided that
the small ball estimates are tight. To obtain tight small ball estimates is in itself a difficult
problem in general. If our estimates of the $\epsilon$-approximations are tight, our entropy bounds are tight modulo log-terms, as is explained in 
\cite{van2023lasso}. 
We conjecture that there are indeed superfluous 
log-terms in our bounds. We will see in Subsection \ref{q0=1.section} however, that we can recover the (at the time
of writing to the best of our knowledge) so far existing
entropy upper bound for the example considered there. 

In our exposition we refrain from the Hilbert/Banach space formulation to avoid too many abstractions.
We discuss here an be extension to
function spaces $\Psi$ in $ L_2 (\mu)$ where 
$\mu$ is a measure on
$\R$.  
Let $\Theta$ be some parameter space. To avoid digressions, suppose that $\Theta $
is Euclidean space or a subset thereof. We let for each $v \in \Theta$
$$ \psi(\cdot , v) = \psi_v  \in L_2 (\mu)  . $$
Define
$$ \Psi:=  \psi (\cdot , \cdot) = \{ \psi (u,v) : \ (u, v) \in \R \times \Theta  \} $$  
%For $V_1 \times V_2 \subset \Theta^2 $ we write
% $$ \Psi_{V_1} \otimes \Psi_{V_2} := \{ \psi (u_1, v_1) \psi (u_2, v_2)  : \ u = (u_1 , u_2)  \in \R^2, 
% (v_1, v_2)  \in V_1 \times V_2  \} . $$
 Let ${\bf V}$ be the collection of all finite subsets of $\Theta$.
Then
${\rm absconv } (\Psi \otimes \Psi ) $ is defined as the closure of the set
$$ \biggl \{ \sum_{v_1 \in V_1 }\sum_{v_2 \in V_2 }  \psi(u_1 , v_1 ) \psi(u_2, v_2) \beta_{v_1, v_2 } : $$
$$ (u_1 , u_2) \in [0,1]^2 , \ 
\beta_{v_1 , v_2 } \in \R, \ 
\sum_{v_1 \in V_1 , v_2 \in V_2 } | \beta_{v_1, v_2 }| \le 1 , (V_1 , V_2) \in {\bf V}^2 \biggr  \} . $$
For $d > 2$, the set
$${\cal F} := {\rm absconv } (\underbrace{\Psi \otimes \cdots , \otimes  \Psi }_{ d \ {\rm times} } ) $$
is defined in an analogous manner, taking $(u_1 , \ldots , u_d) \in \R^d $, and
$(v_1 , \ldots , v_d)$ in $\Theta^d$ or in finite subsets $V_1 \times \cdots \times V_d $ thereof. 

We can now state a function space version of Corollary \ref{dfold.corollary}. Let us first define
the nested sequence of approximations.

\begin{definition} \label{nested2.definition} Let, for some $q \in \Nat$, $\{ {\cal V}_k \}_{k\in \Nat}$ a nested sequence of linear spaces in $L_2(\mu)$
with ${\rm dim} ({\cal V}_k) = q 2^k $, \ $k \in \Nat_0 $. Suppose that for some $\gamma>0$
and $W>0$  and for all $k\in \Nat$,
$$  \max_{v \in [0,1] }  \min_{f \in {\cal V}_k }\int ( \psi_v - f )^2 d\mu \le \gamma^2 2^{-2k/W}. $$
Then $\{ {\cal V}_k \}_{k\in \Nat_0}$ is called a nested sequence of approximations of $\Psi$
with parameters $(q, \gamma , W)$. 
\end{definition}

\begin{theorem} \label{functionentropy.theorem} 
Suppose $\Psi:= \{ \psi_v: \ v \in [0,1] \} \subset L_2 (\mu)$ is a collection of functions with
$L_2(\mu)$-norm at most 1, and with for some constant $V>0$
$$N (\epsilon, \Psi) \lesssim \epsilon^{-V} . $$
Let $\{ {\cal V}_k \}_{k \in \Nat_0 } $ be a nested sequence
of approximations with parameters $(q, \gamma, W)$. Then, with
$${\cal F} := {\rm absconv } (\underbrace{\Psi \otimes \cdots , \otimes  \Psi }_{ d \ {\rm times} } ) \subset
L_2 (\mu \times \cdots \times \mu),
$$
we have
$$ {\cal H} (\epsilon, {\cal F} ) \lesssim \epsilon^{-{2W \over 2+W}} \log^{d-1} (1/\epsilon) \log^{W \over 2+W} (1/\epsilon) .$$
\end{theorem}

As will be clear from the next subsection, the function space formulation can be easier, as
it can be easier
to evaluate the integrals $\int (\psi_v - f)^2 d \mu$ than discrete versions thereof (see also Section \ref{conclusion.section}
for some details on a discrete version).

\section{The truncated power basis and Vitali total variation} \label{piecewise.section}
In this section we examine the situation for the truncated power 
basis, which has the functions 
$$ \psi (u,v) = (u-v)_+^{j-1}  , \ (u,v) \in [0,1]^2, $$
where $j \in [1 : q]$ for some $q \in \Nat$.  
It includes as special case the
 heaviside functions ${\rm l} \{ \cdot \ge v \}$ and the hinge functions $(\cdot-v)_+$.
The basis is related to the falling factorial basis, see \cite{wang2014falling}. Linear combinations
of these functions are polynomials of degree $q-1$ and linear combinations of products of these functions are called multivariate adaptive regression splines (MARS) in
 \cite{friedman1991multivariate}.  
 
 Define now for a given $q \in \Nat$ and for  a function $g: [0,1]^d \rightarrow \R$
 $$ Dg (u) := \prod_{k=1}^d (\partial^{q} g(u))/ (\partial^{q} u_k ), \ u = (u_1 , \ldots , u_d) \in [0,1]^d  $$
 (whenever it exists).
 Then $\| D g \|_1:= \int | Dg | d (\mu \times \cdots \times \mu)$ is called the $q$-th order Vitali total variation of $g$. (In fact, one can
 define $\| D g \|_1$ as 
 the first order Vitali total variation of $ \prod_{k=1}^d (\partial^{q-1} g(u))/ (\partial^{q-1} u_k )$, where the first
 order Vitali  total variation  of a function $h$ is $\int |d h| $.)
 We can write
 $$ g = f_{0,g}  + f_g , $$
 where $f_{0,g} $ is in the null space of $D$: $Df_{0,g} =0$ and $f_g$ is
 orthogonal to this null space. We note that this null space is infinite-dimensional
 as soon as $d>1$. We will neglect the null space of $D$ at first. In view its infinite-dimensionality
 for $d>1$, it will be useful to also assume some regularity for it when it can not be neglected.
 See Subsection \ref{q0=1.section} and Subsection \ref{q0=2.section} for more details in the cases considered there.
 Let $\Psi = \{ \psi_v = (\cdot - v)^{q-1} : \ v \in [0,1]\} $, and
 $${\cal F} = {\rm absconv} ( \Psi \otimes \cdots \otimes \Psi )  . $$
 Then
 $${\cal F}=\{ f_g : \| Df_g \|_1 \le 1 \}  $$
 is the set of all functions orthogonal to the null-space of $D$ with $q$-th order Vitali total variation
 bounded by 1. 
 We take as metric for the approximations of $\Psi$ the metric endowed by the $L_2 (\mu)$-norm,
 where $\mu $ is Lebesgue measure on $[0,1] $. 
 Then indeed the functions $\psi_v$ have $L_2 (\mu)$-norm at most one, and for $q > 1$ the norm is
 strictly smaller than 1. One may renormalize to gain in the parameter $\gamma$ 
 in the nested approximations. However, to avoid carrying around the normalization all the time, we keep the 
 $\psi_v$ as defined above. 
% We define the monomials
% $$b_j (u)=  u^{j-1} , u \in [0,1], \ j\in [1 : q] .$$

 \begin{lemma}\label{nestedMARS.lemma} Let 
 $\Psi:= \{ \psi_v (\cdot) = (\cdot -v)_+^{q-1} : \ v \in [0,1] \} $. Then there is a sequence
 of nested approximations $\{ {\cal V}_k \}_{k \in \Nat_0}$ as described in Definition \ref{nested2.definition},
 with parameters $(q , \gamma , W)$
 where $\gamma =  {1 / \sqrt {2 q-1}}$ and
  $W= 2/ (2 q -1 )$. To be precise, we can take
  $$ {\cal V}_k ={\rm span} \biggl ( \biggl  \{ (\cdot - (l-1)2^{-k} )_+^{j-1} : j=[1: q] , \ l=[1: 2^{k} ] \biggr \} \biggr ),\ k \in \Nat_0.$$
\end{lemma}

Recall that the proofs for this section can be found in Section \ref{proofspiecewise.section}.

\begin{corollary}\label{MARS.corollary}
Let $\Psi:= \{ \psi_v (\cdot) = (\cdot -v)_+^{q-1} : \ v \in [0,1] \} $. 
Then $N (\epsilon , \Psi ) \lesssim 1/ \epsilon $. 
  By Theorem \ref{functionentropy.theorem},
  for ${\cal F} := {\rm absconv} ( \Psi \otimes \cdots \otimes \Psi) $ with
  the metric endowed by the $L_2 (\mu \times \cdots \mu )$-norm for functions on $[0,1]^d$,
  the entropy can be bounded by
  \begin{equation}\label{entropyMARS.equation}
  {\cal H} ( \epsilon , {\cal F} ) \lesssim \epsilon^{-{1 \over q}}  \log^{d-1} (1/\epsilon)\log^{1 \over 2 q} (1/\epsilon).
  \end{equation}
  \end{corollary}
  
  \begin{remark} For $d=1$ the above bound recovers the one as given in \cite{Birman:67} up
  to the log-term $\log^{1 \over 2q} (1/\epsilon)$. 
   \end{remark}

  \subsection{The case $q=1$} \label{q0=1.section}This is the situation of absolute convex combinations
  of heaviside functions. The class ${\cal F}$ is then the class of all
  functions that can be written as the difference of two $d$-dimensional distribution functions, scaled by $1/2$.
  In \cite{blei2007metric} one finds the lower and upper bound
  $$ \epsilon^{-1}  \log^{d-1} (1/\epsilon) \lesssim {\cal H} ( \epsilon , {\cal F} ) \lesssim \epsilon^{-1}  \log^{d-1/2} (1/\epsilon)$$
  and we see that  their upper bound coincides with our upper bound in (\ref{entropyMARS.equation}).
  Let us look at the case $d=2$ in some more detail.
  For $f \in {\cal F}$ one has $f(u_1, 0) = f(0,u_2) = f(0,0)=0$.
 In general for a function $g:[0,1]^2 \rightarrow \R$ it is true that for $u = (u_1 , u_2) \in [0,1]^2$
  \begin{eqnarray*}
   g(u) &= & g(u_1,0) + g(0, u_2) - g(0,0) \\
   &+&  g(u_1, u_2) -g(u_1,0) - g(0, u_2) + g(0,0)   \\
   &:=& g(u_1,0) + g(0, u_2) - g(0,0)  + f(u) ,
  \end{eqnarray*}
  with $f \in {\cal F}$.
  The constant term $g(0,0)$ is typically not restricted (in statistical applications).
  So one may extend the situation to get entropy bounds of the same order for the class\footnote{
  We remark that for $h$ a function on $[0,1]$, the integral
  $\int |h^{\prime} | d \mu$ can be read as $\int |d h | $, that is, as its total variation. Then
  differentiability of $h$ is not required. }
  $$ {\cal G} :=\biggl  \{ g: u \mapsto  g_{1,0}(u_1) + g_{0, 1}  (u_2) + f(u): \ u= (u_1 u_2)  \in [0,1]^2 ,$$ $$
   f \in {\cal F}, \ \int ( |g_1^{\prime} | + | g_2^{\prime} | )d \mu ,
   \le 1   \biggr \} . $$
  One can take the functions in ${\cal G}$ to have mean zero (so that
  they are orthogonal to the constant function) and also one may take all functions involved
  (the constant function, 
  $g_1$, $g_2$ and $f$) mutually orthogonal (an ANOVA decomposition). 
  
  \subsection{The case $q =2$}\label{q0=2.section} 
  In this subsection,  $\psi_v ( u) = (u-v)_+ $, $(u,v) \in [0,1]^2$.
  For $g$ a (twice differentiable) function on $[0,1]$ 
  %$$ g(u) = g(0) + \int_0^u g^{\prime} (x) dx   $$ 
  $$ g(u)= g(0) + u g^{\prime} (0) + \int_0^1 (u-v)_+ g^{\prime \prime} (v) dv , \ u \in [0,1].$$
   In two dimensions this becomes: for $g$ a function on $[0,1]^2$ , and $u= (u_1 , u_2) \in [0 , 1]^2$, and
   whenever the derivatives involved in the expression exist,
%  $$ g(u_1 , u_2) = g(0, u_2) + u_1 g_{U_1}  (0, u_2) ) + \int_0^1 (u_1 - v_1)_+  g_{U_1U_1} (v_1, u_2 ) dv_1 $$
%  $$ g_{U_1U_1} (v_1, u_2 )= g_{U_1U_1} (v_1, 0) + u_2 g_{U_1U_1, U_2} (v_1 , 0) $$
%  $$ + \int_0^1 (u_2 - v_2)_+  g_{U_1U_1, U_2 U)2} (v_1, u_2 ) dv_2 $$
%  Thus
%  $$  \int_0^1 (u_1 - v_1)_+  g_{U_1U_1} (v_1, u_2 ) dv_1 $$
%  $$ = \int_0^1 (u_1 - v_1)_+g_{U_1U_1} (v_1, 0) dv_1 + \int_0^1 (u_1 - v_1)_+ u_2 g_{U_1U_1, U_2} (v_1 , 0) dv_1
%  + \int_0^1 (u_1 - v_1)_+\int_0^1 (u_2 - v_2)_+  g_{U_1U_1, U_2 U_2} (v_1, u_2 ) dv_2 $$
%  $$ = g(u_1 , 0) - g(0,0) - u_1 g_{U_1}  (0,0) $$
%$$ + u_2 g_{U_2} (u_1 , 0) - u_2 g_{U_2} (u_1 , 0) - u_1 u_2 g_{U_1U_2} (0,0) $$
%$$  + \int_0^1 (u_1 - v_1)_+\int_0^1 (u_2 - v_2)_+  g_{U_1U_1, U_2 U_2} (v_1, u_2 ) dv_2 $$
%and so
%$$ g(u_1,u_2) = g(0, u_2)  + u_1 g_{U_1} (0,u_2) + g(u_1 , 0) - g(0,0) - u_1 g_{U_1}  (0,0) $$
%$$ + u_2 g_{U_2} (u_1 , 0) - u_2 g_{U_2} (u_1 , 0) - u_1 u_2 g_{U_1 U_2} (0,0)  $$
%$$+ \int_0^1 (u_1 - v_1)_+\int_0^1 (u_2 - v_2)_+  g_{U_1U_1, U_2 U_2} (v_1, u_2 ) dv_2 $$
\begin{eqnarray*}
 g(u) &= &a_0 + a_1 u_1 + a_2 u_2 + a_{1,2}   u_1 u_2 \\
 &+ &  g_{0,1} (u_2) + g_{1,0}  (u_1)+ u_1 g_{0,2} (u_2) + u_2 g_{2,0} (u_1) \\
 &+& f(u) , \ u= (u_1 , u_2) \in [0 , 1]^2 , 
 \end{eqnarray*}
 where
 $$ f(u) = \int_0^1 \int_0^1 (u_1 - v_1)_+ (u_2 - v_2)_+ { \partial^4 g(v_1 , v_2) \over  \partial v_1^2 \partial v_2^2 }  dv_1 dv_2 . $$
 Let $b_1\equiv 1$ and let $b_2$ be the identity function. Let ${\cal V}_0 := \{ b_1 , b_2 \}$ and $g^{\bot}$ the part of the function $g$ orthogonal to
 ${\cal V}_0 \times {\cal V}_0$. We define
 $${\cal G} := \biggl \{ g^{\bot}  : \ u \mapsto g_{1,0} (u_1) + g_{0,1} (u_2) + u_1 g_{0,2} (u_2) + u_2 g_{2,0} (u_1) 
 + f(u_1,u_2) : $$ $$ u = (u_1,u_2) \in [0,1]^2 , \  f \in {\cal F} , 
 \ \int\biggl (   | g_{1,0}^{\prime \prime} | + | g_{0,1}^{\prime \prime} | +
 | g_{0,2}^{\prime \prime} | + | g_{2,0}^{\prime \prime} | \biggr ) d \mu \le 1\biggr  \} $$
 with ${\cal F} := {\rm absconv} ( \Psi \otimes\Psi) $. Again, one can take all functions
 mutually orthogonal.

\subsection{Multi-resolution analysis} \label{multiresolution.section}
Although not necessary for the entropy results, it is of interest to see that
nested orthogonal bases $\{ p_j \}_{j =1}^{q 2^k } $ for
the nested linear spaces $ \{ V_k \}$ can be constructed as a wavelet. For the case $q=1$ we then get
- not surprisingly - the Haar basis. In general we get an orthogonal basis $P= \{ p_1 , p_2 , \ldots \}$
which can serve as multi-resolution, piecewise polynomial dictionary for the functions in $L_2 (\mu)$ that can be well-approximated by
piecewise polynomials of degree $q-1$.

Define for $F:\ [0,1] \rightarrow \R^{q_1}$ and $G: \ [0,1] \rightarrow \R^{q_2}$,
$$ \langle F , G \rangle := \int F G^T d \mu  \in \R^{q_1 \times q_2 } .  $$
Let for $u \in [0,1]$, $b_j(u) := u^{j-1}$, $j=1 , \ldots , q$. Then the space ${\cal V}_0$ spanned
by $B_0 := \{ b_1 , \ldots , b_{q} \} $ is location and scale invariant in the sense that for all
changes of location $c_1 \in \R$ and scale $c_2 >0$, it holds that
$${\rm span} \biggl  \{  B_0  \biggl ({\cdot  - c_1 \over  c_2}  \biggr ) \biggr  \} = {\cal V}_0 . $$
We define
$$B_{1,1} := B_{1} := B_0 A_{1,1}{\rm l}_{[0,1/2)}+ B_0 A_{1,2} {\rm l}_{[1/2, 1] } $$
where the  matrices $A_{1,1} $ and $A_{1,2} $, both  in $\R^{q \times q }$,
are such that $\langle B_{1} , B_0 \rangle =0 $, i.e, all functions in $B_{1}$ are orthogonal to 
those in $B_0$.  

Set for $l \in [1 : 2^{k-1}]$ and $k \in \Nat$,
$$ B_{k, l} (u) := B_{1} \biggl ( 2^{k-1} (u - (l-1)2^{-(k-1)} ) \biggr ){\rm l}\{  (l-1) 2^{-(k-1)} \le u < l 2^{-(k-1)} )\} . $$

\begin{lemma}\label{orthogonalB.lemma} We have for all $k \in \Nat$ and all $l \in [1 : 2^{k-1}]
$ 
$$ \langle B_0 , B_{k,l} \rangle =0 $$
and for all $(k_1 , k_2) \in \Nat^2$ with $k_1 \not= k_2$ and all $(l_1, l_2)  \in [1 : 2^{k_1-1}]
\times [1: 2^{k_2-1} ]$ 
$$ \langle B_{k_1, l_1 } , B_{k_2,l_2} \rangle =0 .$$
Moreover, ${\cal V}_0 = {\rm span} (B_0) $ and for all $k \in \Nat$
$${\cal V}_k = {\rm span}\biggl  ( B_0, \biggl \{ \{ B_{k^{\prime},l^{\prime}}\} _{l^{\prime}=1}^{2^{k^{\prime}-1} }\biggr  \}_{k^{\prime}=1}^k \biggr ) $$
with ${\cal V}_k $ given in Lemma \ref{nestedMARS.lemma}. 
\end{lemma}

 We now have basis functions
 $$ \biggl \{ B_0 , B_{1,1}, (B_{2,1}, B_{2,2}) , (B_{3,1}, B_{3,2}, B_{3,3} , B_{3,4}) , \ldots \biggr \},$$
 where the block $B_0$ and each block $B_{k,l}$ consist of $q$ linearly independent functions,
 and all blocks are mutually orthogonal.
 We can do a Gram-Schmidt orthogonalization within all these blocks of size $q$  
 to obtain $q$ orthogonal functions within each block, and then normalize these
 to obtain $q$ orthonormal functions in each block.
 Then we get ${\rm span}(B_0) = {\rm span} (\{ p_1 , \ldots , p_{q}\} )$ and for $l=1: 2^{k-1}$ and $k \in \Nat$
 $$ {\rm span} (B_{k,l} )= {\rm span} (\{ p_{(2^{k -1}+l-1) q +1 } , \cdots ,  p_{(2^{k-1}+l)q}  \}).$$
 The functions $\{ p_1 , p_2 , \ldots \}$ form an orthonormal basis of $L_2 (\mu)$ consisting of piecewise polynomials
 of degree $q-1$. .

\subsection{A piecewise linear multi-resolution basis} \label{hinge.section}
For $q=1$ Lemma \ref{orthogonalB.lemma} shows the Haar basis for functions on
$[0,1]$. Let us consider here in some more detail the case
case $q=2$, which is about the hinge functions
$\psi_v(u) := (u-v)_+ , \ (u,v) \in [0,1]^2 $.
With $q=2$, Lemma \ref{orthogonalB.lemma}  constructs an orthonormal basis where the basis functions are piecewise linear. We give their explicit expressions. 
Let $ b_1 \equiv 1$ be the constant function and, for $u \in [0,1] $, $b_2(u)=u$,
$b_3(u)= {\rm l} \{ u \ge 1/2 \} $ and $b_4 (u) = (u-1/2)_+$. 
Let ${\cal V}_0 $ the space spanned by
$\{ b_1, b_2\} $ and ${\cal V}_1$ the space spanned by $\{ b_1 , b_2, b_3, b_4 \} $. 
We perform a Gram-Schmidt
orthogonalization. Let, for $k \in [1 : 3 ]$, $\Pi_k $ be the projection operator onto the
space spanned by $b_1 , \ldots , b_k$. We take $g_1= b_1 $, $g_2 = (b_2- \Pi_1 b_2) $,
$g_3 = (b_3 - \Pi_2 b_3 )$ and $g_4 = (b_4 - \Pi_3 b_4 ) $.
Writing for $f \in L_2(\mu)$,  $ \| f \|_{\mu}^2 := \int f^2 (u)  d u$, 
we obtain after the normalization $p_k := g_k / \| g_k \|_{\mu} $, $k\in [1:4]$, 
an orthonormal basis $\{ p_1 , p_2 \}$ for ${\cal V}_0 $
and $\{ p_1 , p_2 , p_3 , p_4 \}$ is an orthonormal basis ${\cal V}_1$. 
We see that - after some calculations - for $u \in [0,1]$, 
$$ p_1 \equiv 1 ,  \ p_2 (u) = (x-1/2)/\sqrt {1/12} , $$
$$ p_3 (u) = \biggl [ {1 \over 4} (1-6u){\rm l} \{ u < 1/2\} + {1 \over 4} (5-6u) {\rm l} \{ u \ge 1/2 \} \biggr ] / \sqrt {7/128} $$
and
$$ p_4 (u) \propto  \biggl [ \biggl ( - {11 \over 28} + {73 \over 28}u \biggr )  {\rm l} \{ x < 1/2\} + 
\bigg[ \biggl ( - {83 \over 28 } + {101 \over 28} u \biggr ) {\rm l} \{ u \ge  1/2 \} \biggr ] 
 . $$
Note that $p_3$ and $p_4$ are piecewise linear, but that the slopes on the two pieces of $p_3$ are equal,
whereas the slopes on the two pieces of $p_4$ are not. 
We now add 4 more basis functions, by splitting the interval $[0,1] $ into $4$ subintervals
$[0, 1/4), [1/4, 1/2), [1/2,  3/4), [3/4 , 1]$.  We consider for $u \in [0,1]$ the functions 
$$b_5 (u) := {\rm l } \{ u\ge 1/ 4\}  , \  b_6 (u) := (u- 1/4)_+ , $$
$$ b_7(u):= {\rm l} \{ u\ge 3/4\} , \ 
b_8 (u) := (u-3/4)_+ . $$
Let us start with the left interval $[0,1/2] $. Note that
$$b_5 (u) = {\rm l } \{ 1/4 \le u <  1/ 2\} + {\rm l} \{ u \ge 1/2 \} . $$
and
$$ b_6 (u) = (u-1/4)_+  {\rm l } \{ 1/4 \le u <  1/ 2\} + (u-1/4)_+{\rm l} \{ u \ge 1/2 \} $$
The functions $u \mapsto {\rm l} \{ u \ge 1/2 \}$ and $u \mapsto (u-1/4)_+ {\rm l} \{ u \ge 1/2\} $ are clearly in the space spanned by
${\cal V}_1$. So we can remove them by using a linear combination of $\{ p_1 , p_2 , p_3 , p_4 \}$.
Now we only have to look at the parts
$$ b_5^L (u) = {\rm l } \{ 1/4 \le u <  1/ 2\}, \ b_6^L (u) := (x-1/4)_+ {\rm l} \{ 1/4 \le u < 1/2\} , \ u \in [0,1] . $$
On the interval $[0 , 1/2) $ the functions $\{ p_1 , \ldots , p_4 \}$ are affine linear,
that is, restricted to $[0 , 1/2) $ the space spanned by $\{ p_1 , \ldots , p_4 \}{\rm l} [0 , 1/2) $  is spanned by the constant and the
linear function on $[0, 1/2)$. 
We can do the same Gram-Schmidt orthogonalization as before, now
at a resolution level $k=2$, for $b_5^L$ and $b_6^L$.  The same construction can be done 
for the right interval $[1/2, 1] $. We have thus constructed  the space ${\cal V}_2$ spanned by
an orthonormal basis
$\{ p_1 , \ldots , p_8 \}$. This process can be repeated to arrive at an orthogonal basis
$\{ p_1 , p_2 , \ldots \} $ for  the nested sequence 
of Lemma \ref{nestedMARS.lemma}
$${\cal V}_k= {\rm span} \biggl ( \biggl \{ {\rm l} \{ \cdot  \ge (l-1)2^{-k} \} , (\cdot -(l-1) 2^{-k})_+ \biggr \}_{l=1}^{2^k} \biggr ) , \ k\in \Nat_0,  $$
with starting value $q=2$.  Notice that ${\cal V}_k$ has $q 2^k$ basis functions, which equals the sum
$q (1+ \sum_{k^{\prime}=1}^{k} {2^{k^{\prime}-1}} )$ of the number of basis functions per resolution level.
We obtain an orthogonal basis $\{ p_1 , \ldots , p_{q2^k} \}$ for ${\cal V}_k $ and the $q 2^{k-1} $ basis functions
at level $k$ are scaled versions on subintervals of the basis functions at level 1.
The orthogonal basis thus obtained $P:= \{ p_1 , p_2 , \ldots \}$ can be used as expansion for functions in $L_2 (\mu)$.
It is a dictionary for functions that are well-approximated by piecewise linear functions.
Moreover, $P \otimes \cdots \otimes P$ is an orthogonal basis for functions in $L_2 ( \mu \times \cdots \times \mu)$,
and is a dictionary for functions that we well-approximated by products of piecewise linear functions of the coordinates. 

%The sequence  $\{ {\cal V}_k \} $ is a nested sequence of approximations for $\Psi$ with
%parameters $q=2$, $\gamma=1/\sqrt 3$ and $W = 2/3$.  Indeed, for
%$k \in \Nat$, $l\in [1: 2^k]$ and 
%for $v \in [(l-1)2^{-k}, l 2^{-k})$, we have
%$$ \psi_v (u) - \biggl ( (u- l 2^{-k} )_+ + 2^{-k} {\rm l} \{ u \ge l 2^{-k} \}\biggr )  =  (u-v)  {\rm l} \{ v \le u < l 2^{-k} \}  . $$
%and 
%$$ \int_v^{l 2^{-k}}   (x-v)^2 dx = {1 \over 3 } (l 2^{-k} -v)^3  \le {1 \over 3} 2^{-3k} .  $$ 
%
%

\section{Conclusion} \label{conclusion.section}
We derived entropy bounds for absolute convex hulls using low-dimensional approximations 
of the extreme points, and applied this to the case where the extreme points are tensors in $L_2 (Q \times \cdots \times Q)$. The result can be extended to tensors in $L_2 ( {\bf Q})$ where ${\bf Q}$ 
and a  product measure are mutually dominating. We started out with the vector space formulation and then stated a case where
one deals with functions instead of vectors as obvious extension. 

We studied the
truncated power series $\Psi$ as special case. This is where the functions are in $L_2 (\mu \times \cdots \times \mu)$
with $\mu$ Lebesgue measure on $[0,1]$. Let us have a brief look at the discrete setting, 
where instead of Lebesgue measure on $[0,1]$ one has the
uniform measure $Q_m$ on $[1 : m ] $. The squared $L_2 (Q_m)$-norm is then
 $$ \| v \|_{Q_m}^2 = \sum_{i=1}^m v_i^2 / m , \ v \in \R^m.  $$

\vfill\eject
For the case $q=2$, 
the vectors $\{ \psi_j\}_{j=3}^m$ of the matrix $\Psi$ are 
$$ \psi_{i,j} := {1 \over m} (i- (j+1)){\rm l} \{ i \ge j \} , i \  \in [1:m] . $$
Further details are in the footnote\footnote{
Then clearly $\psi_j \in L_2 (Q_m)$ for all $j$.
For $\beta := ( \beta_3 , \ldots , \beta_m)^T \in \R^{m-2}$ and
$ f = \Psi \beta $
we have
$$f_1=f_2=0, \  f_j= ((j-2) \beta_3 + \ldots + \beta_{j})/m , \ j\in [3: m ] $$
and
$\beta_j = m (f_{j} - 2 f_{j-1}+ f_{j-2} ) , \ j\in [3 : m ] $. Thus
$$ \| D f \|_1 = m\sum_{j=3}^m | f_{j} - 2 f_{j-1}+ f_{j-2}   | . $$
We can as in the continuous case approximate $\Psi$ by piecewise constant and piecewise
linear functions with jump, respectively kink, on the grid $\{ l 2^{-k} \}_{l=1}^{2^k} $. Let us
present some details. 
Assume (for convenience) that $m= 2^K$ for some $K \in \Nat$. Then for $k \le K$ and
all $l \in [ 1 : 2^k]$ it holds that $j_k:= m l 2^{-k} \in [1 : m ]  $.
For $j/m \in [(l-1)2^{-k}, l 2^{-k}) $ we get for $i \in [1:m ]$
$$ \psi_{i,j} - \underbrace{\psi_{i, j_k  }}_{{\rm piecewise \ linear} }  - \underbrace{(j_k - j) {\rm l} \{ i \ge j_k \} }_{{\rm piecewise \ constant }} = \psi_{i,j} {\rm l} \{ j \le i < j_k \} =: e_{i,j} $$
and
\begin{eqnarray*}
\| e_j \|_{Q_m}^2 &=&
{1 \over m} \sum_{j \le i < j_k } \psi_{i,j}^2 \\
&  \le &{1 \over m^3} \sum_{i=1}^{m 2^{-k} } i^2\\
 &=& {1 \over m^3} { m2^{-k} ( m 2^{-k} + 1) ( 2m 2^{-k} + 1)  \over 6 } . 
\end{eqnarray*}}, showing that nothing new happens as compared to the continuous case, but that
explicit expressions are somewhat more involved. 

Vitali total variation can be compared with  Kronecker total variation. In the first case $\prod_{k=1}^d \partial / \partial u_k$
whereas in the second it is $\sum_{k=1}^d \partial / \partial u_k $. For higher orders we get
$D_V= \prod_{k=1}^d \partial^{q}/ (\partial^{q} u_k) $ for the $q$-th order Vitali total variation, and
we can compare this with 
$ D_M:= ( \sum_{k=1}^d \partial  / (\partial  u_k ))^{qd}  $. In the last case we have all
mixed partial derivatives of order $qd $, whereas $q$-th order Vitali total  variation only considers
one of these terms. However, the null space  of $D_M$ is finite-dimensional, whereas for
$d>1$ that of $D_V$ is infinite-dimensional. This is why we need to restrict (or regularize) terms in the null space of $D_V$.
To avoid that the combinatorics muffles the main point, let us look only at the case $q=d=2$.
Then $q d = 4$ and there are 5 mixed partial derivatives of order 4: 
$$\partial^4 / \partial^4 u_1, \  
\partial^4 / (\partial^3u_1 \partial u_2) , \ \partial^4 / (\partial^2 u_1\partial^2 u_2)  , \ \partial^4 / (\partial u_1\partial^3 u_2) , \partial^4 / \partial^4  u_2 . $$
Only the term $\partial^4 / (\partial^2 u_1\partial^2 u_2)$ is taken into account by second order
Vitali total variation $D_V$, so it misses four terms. However, as we have seen in Subsection \ref{q0=2.section},
we need to restrict four  terms of the infinite-dimensional null space, i.e., the number of terms
to restrict in the null-space is equal to the number of mixed derivatives $D_V$ missed as compared to
$D_M$. In other words, 
although not being equal, $D_M$ and  $\{ D_V+ {\rm null \ space\ restrictions}\}$ are comparable in terms of complexity.

The approach to entropy bounds using low-dimensional approximations is closely related to
results in literature, as explained in Section \ref{functionspaces.section}. This means for instance that
our approximation in Section \ref{piecewise.section} , say for the case $q=2$, can be applied to
get small ball estimates for the $d$-dimensional integrated Brownian sheet. See 
\cite{bass1988probability} who uses the Haar basis for obtaining small ball estimates for the case $q=1$, i.e.,
for the Brownian sheet itself. 

\section{Proof of Theorem \ref{entropy.theorem}.} \label{proofentropy.section}
We first show that
$$ {\cal H}^{\rm loc} ( \epsilon ) \lesssim \epsilon^{-{2W \over 2+W} } \log^{2w \over 2+W} (1/\epsilon) \log^{W \over 2+W} (1/ \epsilon) 
. $$
For convenience, we replace $\epsilon $ by $4 \epsilon$. We take $u\ge \epsilon >0$
to be chosen later.  Set ${\cal V}_u$ as the linear space satisfying
$\delta ({\cal V}_u , X) = M(u, X) $ where  $M(u, X) \le B u^{-W} \log^w (1/u) $ for some constant $B$.
Let $f_0 \in \R^n$ be arbitrary and ${\cal W}_u := {\rm span} ({\cal V}_u , f_0) $. 
Define
$$ X_u^{\bot} := \{ x_1^{{\cal W}_{u}^{\bot } },   \ldots ,  x_p^{{\cal W}_{u}^{\bot } }\} .$$
For any $\beta  \in \R^p$
$$X \beta =  ( X - X_u^{\bot} )  \beta+  X_u^{\bot} \beta . $$
Moreover, for $f = X\beta $ 
$$\| (X - X_u^{\bot} ) \beta - f_0 \|_2 = \| ( X - X_u^{\bot} )  (\beta - f_0 ) \|_2 \le \| f - f_0 \|_2 .$$ 
Thus
\begin{eqnarray*}
& & {\cal H} (4\epsilon, \{f \in {\cal F}: \ \| f - f_0 \|_2 \le 16 \epsilon\} )\\
 &\le & {\cal H} (  \epsilon, {\rm absconv} (X- X_u^{\bot}) ) + {\cal H}( 3 \epsilon , {\rm absconv} (X_u^{\bot} ) )\\
&\le& (M(u,X) +1)  \log ( 3 \times 16 ) + {\cal H}( 3 \epsilon , {\rm absconv} (X_u^{\bot} ) ),
\end{eqnarray*}
where in the second inequality we used 
$${\rm absconv} (X- X_u^{\bot}) \subset
\{ f \in {\cal W}_u : \ \| f - f_0 \|_2 \le 16  \epsilon \} , $$
with ${\cal W}_u$ a linear space with ${\rm dim} ( {\cal W}_u ) \le M (u , X) +1 $,  and where we applied Lemma 14.17 in \cite{BvdG2011}.

To bound the entropy of ${\rm absconv} (X_u^{\bot})$ we first replace the columns of $X_u^{\bot}$ by
a subset with size not depending on $p$ but on $\epsilon$.
Let $\{ \phi_k \}_{k=1}^N$ be a minimal
$\epsilon$-covering set of $X$, where $N= N( \epsilon, X)$ with 
$N(\epsilon , X) \le A \epsilon^{-V} $ for a constant $A$. Then
$ \{ \phi_k^{{\cal W}_u^{\bot} }\}_{k=1}^N$ is clearly an $\epsilon$-covering set
of $X_u^{\bot}$. We can replace this by a $(2 \epsilon)$-covering set
$X_{u, \epsilon}^{\bot} := \{ x_{j_1}^{{\cal W}_u^{\bot}} , \ldots , x_{j_N}^{{\cal W}_u^{\bot}} \} \subset X_u^{\bot} $, so that
$\max_{1 \le k \le N} \| x_{j_k}^{{\cal W}_u^{\bot}} \|_2 \le u $. 
We get
$${\cal H}( 3\epsilon , {\rm absconv} (X_u^{\bot} ) )\le {\cal H} ( \epsilon , {\rm absconv} (X_{u,\epsilon} ^{\bot} ) ). $$

 By Maurey's Lemma (see \cite{carl1985inequalities} or Lemma  9 in \cite{van2023lasso}) we know that for $0 < \epsilon \le u$  the set ${\rm absconv} (X_{u,\epsilon}^{\bot} )$ has
$\epsilon$-entropy
$$ {\cal H} (\epsilon , {\rm absconv} (X_{u,\epsilon}^{\bot} )) \le \biggl \lceil { u^2 \over \epsilon^2 }  \biggl \rceil \biggl (1+ \log (2N+1)\biggr ) .$$
As for the choice of $u$, we trade off the bound for 
$ M( u , X) $ which is a decreasing function of $u$,  against the increasing function $u \mapsto u^2 /\epsilon^2$.
In fact we take
$$u= \epsilon^{ 2 \over 2+W} \log^{w \over 2+W}(1/\epsilon) \log^{-{1 \over 2+W}} (1/\epsilon) 
  . $$
Then we have
\begin{eqnarray*}
& & {\cal H} (4\epsilon, \{f \in {\cal F}: 
\ \| f - f_0 \|_2 \le 16  \epsilon\} )\\
 &\le& \biggl ( B u^{-W} \log^w (1/u) +1\biggr )  \log (3\times 16) 
 \\ &+& 
 \biggl \lceil { u^2 \over \epsilon^2 }  \biggl \rceil \biggl (1+ \log (2A2^V\epsilon^{-V} +1) \biggr ) \\
 & \le & (B+V) \epsilon ^{ - {2W \over 2+W}} \log^{2w\over 2+W}(1/\epsilon) \log^{W \over 2+W}  (1/ \epsilon)
  \\
 & + & {\rm equal \ or \ smaller \ order \ terms}.
 \end{eqnarray*} 
 Since $f_0$ is arbitrary, we can conclude that
 \begin{eqnarray*}
 {\cal H}^{\rm loc} (4\epsilon, {\cal F}) &\le& (B+V) \epsilon ^{ - {2W \over 2+W}} \log^{2w\over 2+W}(1/\epsilon) \log^{W \over 2+W}  (1/ \epsilon)
  \\
 & + & {\rm equal \ or \ smaller \ order \ terms}.
 \end{eqnarray*}
 and hence
  $$ {\cal H}^{\rm loc} (\epsilon, {\cal F} ) \lesssim \epsilon ^{ - {2W \over 2+W}} \log^{2w\over 2+W}(1/\epsilon) \log^{W \over 2+W}  (1/ \epsilon) . $$
  The result for the global entropy ${\cal H}(\epsilon, {\cal F} )$ then follows from (\ref{localentropy.equation}) combined
  with Lemma \ref{localentropy.lemma} below. 

\hfill $\sqcup \mkern -12mu \sqcap$ 

The lemma below allows one to go from local entropies to global ones using the result 
(\ref{localentropy.equation}) of \cite{yang1999information}.

\begin{lemma} \label{localentropy.lemma}
Let $h^{\rm loc} : \ \epsilon \mapsto h^{\rm loc} (\epsilon )$ be the map
$$ h^{\rm loc} (\epsilon ) = B \epsilon^{- A} \log^{\alpha} (1/ \epsilon)$$
where $B>0$, $A >0$ and $\alpha \ge 0 $. Suppose $h$ satisfies $h(1)\le 1$ and for all $k \in \Nat$
$$ h( 2^{-k} / 4) \le h(2^{-k} ) + h^{\rm loc} ( 2^{-k} ) $$
Then for all $K \in \Nat$
$$h(2^{-2K} /4) \le B  \biggl (  2^{2(K+1) A } - 1  \biggr ) {  \log^{\alpha} (2^{2K} )\over 2^{2A} -1 }
. $$
\end{lemma}

{\bf Proof of Lemma \ref{localentropy.lemma}.}
We have
\begin{eqnarray*}
h (2^{-2K} / 4   ) &\le & h( 2^{-2K} ) + B 2^{2KA}  \log^{\alpha} (2^{2K} )\\
&\le& h( 2^{-(2K-2}) + B 2^{(2K-2)A}  \log^{\alpha} (2^{2(K-2)})+
B 2^{2KA} \log^{\alpha} (2^{2K}) \\
& \le & B\sum_{k=0}^{K}  2^{2(K-k) A }  \log^{\alpha} 2^{2(K-k)} \\
& \le & B   \log^{\alpha} 2^{2K} \sum_{k=0}^K 2^{2(K-k) A }  \\
& =& B  \biggl ( 2^{2(K+1) A } - 1  \biggr )  {\log^{\alpha} (2^{2K} ) \over 2^{2A} -1}
\end{eqnarray*}
\hfill $\sqcup \mkern -12mu \sqcap$

\section{Proof of Theorem \ref{tensorcovering.theorem}.}\label{proof.section}

We start with two elementary lemmas.

\begin{lemma} \label{lessthanK.lemma}Let for some $K \in \Nat$, $\{ c_k \}_{k=0}^K$ and $\{ d_k \}_{k=0}^K $ be positive sequences with, for some $q >0$,
$\sum_{l=0}^k c_l \le q 2^k$ and $\sum_{l=0}^k d_l  \le q 2^k $, $k \in [0 : K ] $. Then
$$ \sum_{k_1 + k_2 \le K} c_{k_1} d_{k_2} \le q^2 (K+1) 2^K . $$
\end{lemma}

{\bf Proof of Lemma \ref{lessthanK.lemma}.} We first observe that 
$c_k \le \sum_{ l =0}^k c_{ l} \le q 2^k $, $ k\in [0:K] $. 
Then
\begin{eqnarray*}
 \sum_{k_1 + k_2 \le K } c_{k_1} d_{k_2} &=&
c_0 ( d_1 + \cdots +d_{K} ) + c_1 (d_0 + \cdots  + d_{K-1}) + \cdots + c_{K} d_0\\
&\le& q  2^0 \times q 2^K + q 2^1 \times q 2^{K-1} + \cdots  + q 2^K \times q 2^0\\
&=& \sum_{k=0}^K q^2 2^K \\
&=& q^2 (K+1) 2^K . 
\end{eqnarray*}
\hfill $\sqcup \mkern -12mu \sqcap$

\begin{lemma} \label{greaterthanK.lemma}
Let $\{ c_k \}_{k\in \Nat_0}$ and $\{ d_k \}_{k\in \Nat_0} $ be positive sequences satisfying
$$ \sum_{k \ge 0 } c_k \le 1 , \ \sum_{k \ge 0 } d_k \le 1, $$
and, for some $\gamma > 0$, $W > 0$,
$$\sum_{l > k}  c_l \le \gamma^2 2^{-2k/W} , \ \sum_{l>k}  d_l \le \gamma^2 2^{-2k/W} , \ k \in \Nat_0 . $$
 Then for all $K \in \Nat_0$
$$ \sum_{k_1 + k_2 > K} c_{k_1} d_{k_2} \le 2 \gamma^2 2^{-2K/W} + \gamma^4 K  2^{-2(K-1)/W} . $$

\end{lemma}

{\bf Proof of Lemma \ref{greaterthanK.lemma}.}
Note that $c_0 \le \sum_{k\ge 0 } c_k \le 1$ and for all $k \in \Nat_0$,
$ c_{k+1} \le \sum_{l>k} c_l \le \gamma^2 2^{-2k/W } $. 
We write
\begin{eqnarray*}
\sum_{k_1 + k_2 > K} c_{k_1} d_{k_2} &=& (\sum_{k \ge K+1} c_{k} )(\sum_{k\ge 0 } d_k ) +c_0 \sum_{k \ge K+1} d_k\\
&+&  c_1 \sum_{k\ge K} d_k + c_2 \sum_{k \ge K-1} d_k + \cdots +  c_K \sum_{k\ge 1 } d_k \\
& \le &  \gamma^2 2^{-{2K/W}}  \times 1 + 1 \times \gamma^2 2^{-2K/W } \\
&+&   \gamma^2 2^{-2/W } \gamma^2 2^{-2(K-2)/W } +  \cdots + 
 \gamma^2 2^{-2(K-1)/W} \times \gamma^2 \\
&=& 2 \gamma^2 2^{-2K/W} + \gamma^4 K2^{-2(K-1)/W }.\\
\end{eqnarray*}
\hfill $\sqcup \mkern -12mu \sqcap$

\begin{lemma} \label{projections.lemma} Suppose there exists a nested sequence of approximations for $\Psi$, with parameters
$(q,\gamma,W)$.  Then for all $K \in \Nat$ there is a space ${\cal V}$ of dimension at most
$(K+1) q^2 2^K$ such that
$$ \delta ( {\cal V} , \Psi \otimes \Psi) \le \sqrt {2\gamma^2 2^{-2K/W} +   \gamma^4 K 2^{-2(K-1) /W} }. $$
\end{lemma}

{\bf Proof of Lemma \ref{projections.lemma}.}
Let $k \in \Nat_0$ and ${\cal V}_k = ( p_1 , \ldots , p_{q 2^k})$ be an orthonormal basis of ${\cal V}_k$. Then 
$$P := ( \underbrace{\underbrace{ \underbrace{p_1 , \ldots , p_{q}}_{{\cal V}_0} , p_{q+1}  , 
\ldots p_{2 q} }_{{\cal V}_1}, \ldots   ,p_{q 2^{k-1}+1} ,  \ldots , p_{q2^k} }_{{\cal V}_k} , \ldots  \}.$$
is an orthonormal basis of ${\rm span} (\Psi)$. 
By our assumption on the nested sequence of approximations, it holds for
each $j $, that the coefficients
$ a_j := (a_{1,j} , a_{2,j} , \ldots )^T := P^T \psi_j $ of the vector $\psi_j$ satisfy the tail bound
$$\sum_{l \ge  q2^{k}+1}  a_{l,j}^2 \le \gamma^2 2^{-2k/W} , k \in \Nat. $$

The matrix $ P \otimes P$ is an orthonormal basis of ${\rm span} ( \Psi \otimes \Psi)$:
$$ (P \otimes P)^T ( P \otimes P) = (P^T \otimes P^T ) ( P \otimes P) = (P^T P) \otimes P^T P =
I \otimes I . $$
Let $\psi_{j_1} \otimes \psi_{j_2} $ be one of the columns of $\Psi \otimes \Psi$.
Then for  $a_{j_1} :=  P^T \psi_{j_1}$, $a_{j_2} :=  P^T \psi_{j_2}$
$$ (P \otimes P)^T(\psi_{j_1} \otimes \psi_{j_2} ) = P^T \psi_{j_1} \otimes P^T \psi_{j_2} := a_{j_1} \otimes a_{j_2}  $$
and
\begin{eqnarray*}
 \| a_{j_1} \otimes a_{j_2} \|_2^2 &=& \sum_{l_1\ge 0 } \sum_{l_2\ge 0 } a_{l_1 , j_1}^2 a_{l_2, j_2}^2 \\
 &=&
(\sum_{l_1\ge 0 } a_{l_1, j_j}^2 ) (\sum_{l_2\ge 0 } a_{l_2, j_2}^2 ) \\
&=&
( b_{0, j_1}^2 + b_{1,j_1} ^2 + \cdots)( b_{0, j_2}^2 + b_{1, j_2}^2 + \cdots ) \\
& = & ( \sum_{k_1 \ge 0 } b_{k_1, j_1) }^2 ) ( \sum_{k_2 \ge 0 }  b_{k_2, j_2}^2 )\\
&=& \sum_{k_1\ge 0 }\sum_{k_2 \ge 0 } b_{k_1, j_1 }^2 b_{k_2, j_2 }^2
\end{eqnarray*}
where, for $j \in \{ j_1 , j_2 \}$,
$$ b_{0,j}^2 = \sum_{l=1}^{q} a_{l, j}^2 ,  \ b_{k,j}^2 =\sum_{l=q 2^{k-1} +1 }^{q 2^k} a_{l, j}^2 , k \in \Nat  . $$
Note that 
$$ \sum_{\tilde k  > k} b_{\tilde k, j }^2 = \sum_{l\ge  q 2^{k} +1 } a_{l, j}^2 \le \gamma^2 2^{-2k / W} , \ k \in \Nat_0.$$

It follows that from Lemma \ref{lessthanK.lemma} that
$$ \{ p_{k_1} \otimes p_{k_2} \}_{k_1 + k_2 \le K } $$ has at most
$q^2 (K+1) 2^K $ basis functions. Moreover, by Lemma \ref{greaterthanK.lemma}
$$ \sum_{k_1 + k_2 > K } b_{k_1 , j_1}^2 b_{k_2 , j_2}^2 \le 2\gamma^2 2^{-2K/W} +   \gamma^4 K2^{-2(K-1) /W} . $$
\hfill $\sqcup \mkern -12mu \sqcap$

{\bf Finishing the proof of Theorem \ref{tensorcovering.theorem}.}
Take $K_{\epsilon} \in \Nat$ as the smallest value of $K$ satisfying 
$$2 \gamma^2 2^{-2K/W} + \gamma^4 K 2^{-2(K-1)/W} \le \epsilon^2 . $$
Then we have
$$2^{K_{\epsilon}} \asymp  \epsilon^{-W} \log^{W \over 2} (1/\epsilon) , \ K_{\epsilon} \asymp \log (1/\epsilon) . $$
It follows that 
$$M( \epsilon , \Psi \otimes \Psi) \le (K_{\epsilon}+1) q^2 2^{K_{\epsilon} }\lesssim \epsilon^{-W}  \log^{2+ W\over 2}  (1/ \epsilon) . $$
\hfill $\sqcup \mkern -12mu \sqcap$

\begin{remark} \label{d>2.remark} For $d\ge 2$ we have to find the smallest value of $K_{\epsilon}\in \Nat $
such that $K_{\epsilon}^{d-1} 2^{-{ 2K / W} } \lesssim \epsilon^2 $.  Then we have
$$ 2^{K_{\epsilon} } \asymp \epsilon^{-W} \log^{W(d-1)\over 2} (1/\epsilon) , \ K_{\epsilon} 
\asymp \log (1/\epsilon) . $$
We find for the $d$-fold tensor $X= \Psi \otimes \cdots \otimes \Psi$, 
$$M(\epsilon, X) \lesssim K_{\epsilon}^{d-1} 2^{K_{\epsilon}} \asymp \epsilon^{-W} \log^{(d-1) (2+W) \over 2} (1/\epsilon). $$

\end{remark}

\section{Proofs for the results of Section \ref{piecewise.section}}\label{proofspiecewise.section}
{\bf Proof of Lemma \ref{nestedMARS.lemma}.} 
Clearly 
$$ {\cal V}_k ={\rm span} \biggl (  \{ (\cdot - (l-1) 2^{-k} )_+^{j-1} : j=[1: q] , \ l=[1: 2^{k} ] \} \biggr ),\ k \in \Nat_0$$
is a nested sequence, and ${\cal V}_k$ has cardinality $q 2^k $, $k \in \Nat$, 
% For a given $q$ one can construct the nested subspaces  $\{ {\cal V}_k \}_{k\in \Nat_0} $ as follows.
% Let ${\cal V}_0 :={\rm span} \{ b_1 , \ldots , b_{q} \} $ where we recall that
% $$ b_j (u) = u^{j-1}, 
%  0 \le u \le 1 , \ j\in [1 : q ] . $$
%  Let 
%  $${\cal V}_1 \slash {\cal V}_0 := {\rm span} \{ b_j {\rm l}_{[1/2 , 1 ] }\}_{j=1}^{q} \} , $$
%  $${\cal V}_2 \slash {\cal V}_1 := {\rm span} ( \{ b_j {\rm l}_{[1/4, 1/2)} \}_{j=1}^{q} \} \cup \{ b_j {\rm l}_{[3/4 , 1 ]} \}_{jk-1=1}^{q} \})
%  $$
%  and in general, for $k \in \Nat$
%  $$ {\cal V}_k \slash {\cal V}_{k-1} = \cup_{l \le 2^{-k} \ {\rm even} } {\rm span} ( \{ b_j {\rm l}_{[(l-1) 2^{-k}, l2^{-k} )} \}_{j=1}^{q} \} . $$
%  Clearly this is a nested sequence, and since ${\cal V}_k\slash {\cal V}_{k-1}$ has cardinality $| {\cal V}_k\slash 
%  {\cal V}_{k-1} |= q2^{k-1}$,
%  the cardinality $|{\cal V}_k |$ of ${\cal V}_k$ is $q 2^k$
  as required in Definition \ref{nested2.definition}.   
  As for the approximation, note that for $v \in [(l-2) 2^{-k}, (l-1) 2^k] $, and $v_{l,k} := (l-1) 2^{-k}$
  $$\biggl (\psi_v(u)  - \psi_{v_{l,k}} (u) \biggr ){\rm l}\{ v\le u <(l-1) 2^{-k}\} = (u-v)^{q-1}{\rm l}\{  v\le u <  (l-1) 2^{-k}\}  $$
  and 
  $$ \biggl (\psi_v(u)  - \psi_{v_{l,k}} \biggr  ){\rm l}\{ u \ge (l-1) 2^{-k} \}  $$ $$= \biggl ( (u-v)^{q-1} -
  (u- (l-1) 2^{-k} )^{q-1} \biggr)  {\rm l}\{ u \ge (l-1) 2^{-k} \}  .$$
  The latter shows that
  $$\biggl ( \psi_v  - \psi_{v_{l,k} }\biggr ){\rm l}_{[ (l-1)2^{-k} , 1]}  \in {\rm span}\biggl  (\{ (\cdot - (l-1) 2^{-k} )_+^{j-1} : j=[1: q]  \} \biggr ) .$$
  Since ${\cal V}_k $ contains all the functions $\{ (\cdot - (l-1) 2^{-k})_+^{j-1}   : j=[1: q] , \ l=[1: 2^{k} ] \}$
  and since
  $$ \int_0^{2^{-k}}  u^{2(q-1) }du = {1 \over 2 q-1} 2^{-(2 q-1) k} $$
  we see we have a nested sequence with parameters $q$, $\gamma =  {1 / \sqrt {2 q-1}}$ and
  $W= 2/ (2 q -1 )$.
  \hfill $\sqcup \mkern -12mu \sqcap$

  {\bf Proof of Lemma \ref{orthogonalB.lemma}.} 
Note that for $D_{1,1} (u) := B_1 (u/2)$ and $D_{1,2} (u) := B_1 ((u+1)/2)$,
$$  D_{1,1} (u) = B_0 (u/2) A_{1,1} , \ D_{1,2}  (u) = B_0 ((u+1)/2) A_{1,2} .$$
Hence 
$${\rm span} (D_{1,1} )=  {\rm span} (D_{1,2}) = {\rm span} (B_0).$$
It follows that
$$ \langle B_{1} , D_{1,1} \rangle = \langle B_1, D_{1,2}  \rangle =0  . $$
By definition
$$B_{2,1} (u) = B_{1} (2u) {\rm l}_{[0,1/2)} , \ B_{2,2} (u) := B_{1} (2u-1) {\rm l}_{[ 1/2 ,1 ] }. $$
Therefore
\begin{eqnarray*}
 \langle B_{2,1}  , B_{1} \rangle &=& \int_0^{1/2} B_{1} (2u)  , B_{1}^T(u) du
= {1 \over 2} \int_0^1 B_{1} (u) B_{1}( u/2) dx \\
&=&{1 \over 2} \int_0^1 B_{1}(u) D_{1,1}^T (u) du ={1 \over 2}  \langle B_{1} , D_{1,1} \rangle =0  
\end{eqnarray*} 
and similarly 
\begin{eqnarray*}
 \langle B_{2,2}  , B_{1} \rangle &=& \int_{1/2}^1 B_{1} (2u-1) B_{1}^T (u) du 
 ={1 \over 2}  \int_0^1 B_{1} (u) D_{1,2}^T (u) du \\
 &=&{1 \over 2}  \langle B_{1}  , D_{1,2} \rangle =0 .
 \end{eqnarray*}
Further, for
$$ C_{1,1} (x) := B_0 (u/2), \ C_{1,2} := B_0 ((u+1)/2)  $$
we have
\begin{eqnarray*}
\langle B_{2,1}  , B_0 \rangle &=& \int_0^{1/2} B_{1} (2u) B_0^T (u) du 
 = {1 \over 2} \int_0^1 B_{1} (u) B_0^T (u/2) dx\\
 & =& {1 \over 2} \langle B_{1} , C_{1,2} \rangle=0  
  \end{eqnarray*}
and also
\begin{eqnarray*}
\langle B_{2,2}  , B_0 \rangle &=& \int_{1/2}^1 B_{1} (2u-1) B_0^T (u) dx
 = {1 \over 2} \int_0^1 B_{1} (u) B_0^T ( (u+1)/2) du \\
 &=& {1 \over 2} \langle B_{1} , C_{1,2} \rangle = 0 .
  \end{eqnarray*}
 Repeating the argument we obtain the orthogonality.
 \hfill $\sqcup \mkern -12mu \sqcap$

\bibliographystyle{plainnat}
\bibliography{reference}

\end{document}